\documentclass[12pt]{article}

\usepackage{amsmath,amsthm}
\usepackage{amssymb}
\usepackage{amsfonts}
\usepackage{amscd}

\newtheorem{thm}{Theorem}[section]
\newtheorem{theorem}[thm]{Theorem}

\newtheorem{corollary}[thm]{Corollary}
\newtheorem{lemma}[thm]{Lemma}
\newtheorem{proposition}[thm]{Proposition}

\theoremstyle{definition}
\newtheorem{definition}[thm]{Definition}
\newtheorem{example}[thm]{Example}

\newtheorem{remark}[thm]{Remark}
\newtheorem{observation}[thm]{Observation}

\newcommand{\Z} {\mathbf{Z}}
\newcommand{\C} {\mathbf{C}}

\newcommand{\bm} { {m}}

\newcommand{\salg} {\text{(salg)}}

\newcommand{\sets}{{(\hbox{sets})}} 
\newcommand{\spec}{{\hbox{Spec}}}

\newcommand{\uspec}{\text{\underline{Spec}}}

\newcommand{\Hom}{\mathrm{Hom}}
\newcommand{\rk}{\mathrm {rk}}
\newcommand{\Jac}{\mathrm {Jac}}
\newcommand{\Dim}{\mathrm {dim}}
\newcommand{\rGL}{\mathrm {GL}}
\newcommand{\rSL}{\mathrm {SL}}
\newcommand{\rSp}{\mathrm {Sp}}
\newcommand{\rodet}{\mathrm {odet}}

\newcommand{\rOsp}{\mathrm {Osp}}
\newcommand{\Ber}{\mathrm {Ber}}
\newcommand{\rLie}{\mathrm{Lie}}

\newcommand{\lra} {\longrightarrow}

\newcommand{\del} {\partial}
\newcommand{\defi} {\text def}
\newcommand{\Span}{\hbox{span}}

\newcommand{\ep}{\epsilon}
\newcommand{\ttau}{{\tilde \tau}}

\newcommand{\cO}{\mathcal{O}}
\newcommand{\cH}{\mathcal{H}}

\newcommand{\cK}{\mathcal{K}}
\newcommand{\cB}{\mathcal{B}}

\begin{document}


\bigskip

\centerline{\LARGE \bf  Smoothness of Algebraic} 

\medskip

\centerline{\LARGE \bf Supervarieties and Supergroups }

\bigskip

\centerline{ R. Fioresi}

\smallskip

\centerline{\it Dipartimento di Matematica, Universit\`a di
Bologna }
 \centerline{\it Piazza di Porta S. Donato, 5.}
 \centerline{\it 40126 Bologna. Italy.}
\centerline{{\footnotesize e-mail: fioresi@dm.unibo.it}}

\medskip




\begin{abstract}
In this paper we discuss the notion of smoothness in
complex algebraic supergeometry and we prove that 
all affine complex algebraic supergroups  
are smooth. We then prove the stabilizer theorem in the
algebraic context, providing some useful applications.
\end{abstract}

\bigskip

\section{Introduction}

The category of differentiable supermanifolds was
introduced and discussed in several works among which
\cite{ba, be, ko, le, ma1} from different point of views, 
especially in connection with the important physical
applications, which stem from string theory and ultimately
are related with the problem of the classification of
elementary particles. 

\medskip

In this paper we are interested in algebraic supergeometry
and its relation with its differential counterpart. 
In his foundational work \cite{ma1} on supermanifolds,
Manin defined the notion of superscheme and discussed some
important examples. 

\medskip

Along the same lines we want to understand the 
concept of smoothness in
complex algebraic supergeometry. 
Given the algebraic nature of the problems in the theory
of supermanifolds, we believe that  a deep analysis of
the superalgebraic category can shed light also
on the differential one. Moreover it is the correct category 
to work with, when one wants to discuss quantum deformations 
of the geometric objects.

\medskip

In ordinary algebraic geometry smoothness is a local
notion, strongly linked to the  dimension
of the local ring of the variety at the point. 
Unfortunately, due to the presence of the odd nilpotents,
it is not easy to generalize the idea of
dimension of a ring to the super context. To overcome
this problem, we define smoothness as a property of
the completion of the local superring of the supervariety
at a given point; namely we say that a point
is \textit{smooth} if the local super ring is isomorphic to a power
series super ring. We are then able to show that any 
supervariety admits a unique
supermanifold structure in a
neighbourhood of a smooth point, as in the
classical case, through the application of the implicit
function theorem, after reduction to local complete intersection.

\medskip

Using Cartier's Theorem adapted to supergeometry we
can then prove that all algebraic supergroups are smooth,
in other words, \textit{all affine
algebraic supergroups are also Lie supergroups}.
We apply this result to the case of the stabilizer supergroup
functor of an action of an affine supergroup on an affine supervariety.
After showing that the stabilizer is representable, that
is it is a supergroup, we show that 
classical supergroups are smooth (for a list of classical
supergroups see for example \cite{dm} pg 70).
This fact is generally known,
it is treated for example in a different context
by Gruson in \cite{gr} and by Varadarajan in \cite{vsv} pg 289. 
We however provide 
an independent proof using algebraic tecniques, which we believe
can be of help also in other differentiable supermanifold questions
and can also give other examples of algebraic Lie supergroups.

\bigskip

This paper is organized as follows. 

\medskip

In Section 2, we review
some basic facts of algebraic and differential supergeometry,
among which the definition of 
supermanifolds, supervarieties and their 
functor of points.

\medskip

In Section 3 we give the definition of smooth point of
a supervariety. 
We then prove the super version of the classical result
which states that a smooth point of
a complex algebraic variety admits a supermanifold structure 
in a suitable neighbourhood. 

\medskip

In Section 4 we prove that all (closed) points of
complex algebraic groups are smooth. 

\medskip

In Section 5 we prove the Stabilizer Theorem, which states
that the stabilizer functor for the action of an affine algebraic
supergroup on an affine supervariety is representable by a supergroup
hence it is a smooth variety i.e. a supermanifold.

\medskip

As an application, in Section 6, we show that the classical supergroup
functors as described in \cite{dm} pg 70 are representable, 
i. e. they are algebraic supergroups,
and consequently, they are Lie supergroups.

\bigskip

{\bf Acknoledgements.} We wish to thank Prof. V. S. Varadarajan,
Dr. L. Caston, Prof. D. Gieseker and Prof. M. 
Duflo for helpful comments.







\section{Basic definitions of Supergeometry}

In this section we want to recall some basic definitions
and facts in supergeometry. For more details see
\cite{vsv, cf, dm, ma1}.

\medskip
Let $k$ be the ground field.

\medskip

A {\it superalgebra} $A$ is a $\Z_2$-graded algebra,
$A=A_0 \oplus A_1$, $p(x)$ 
denotes the parity of an homogeneous element $x$.
$A$ is said to be {\it commutative} if
$xy=(-1)^{p(x)p(y)}yx$. $I^{odd}$ denotes the ideal
generated by the odd nilpotents. 

\begin{definition}
A {\it superspace} $S=(|S|, \cO_S)$ is a topological space
$|S|$ endowed with a sheaf of superalgebras $\cO_S$ such
that the stalk $\cO_{S,x}$ is a local superalgebra for all
$x \in |S|$. 
A {\it morphism} $\phi:S \lra T$ of superspaces is given by
$\phi=(|\phi|, \phi^*)$, where $\phi: |S| \lra |T|$ is a 
map of topological spaces and
$\phi^*:\cO_T \lra \phi^*\cO_S$ is a sheaf morphism
such that $\phi_x^*(\bm_{|\phi|(x)})=\bm_x$ where $\bm_{|\phi|(x)}$ 
and $\bm_{x}$ are the maximal ideals in the stalks 
$\cO_{T,|\phi|(x)}$ and $\cO_{S,x}$ respectively.
\end{definition}

The most important examples of superspaces are given by
supermanifolds and superschemes. 

\begin{definition}
Let's consider the superspace $\C^{p|q}=(\C^p, \cH_{\C^{p|q}})$, where
$$
\cH_{\C^{p|q}}|_U=\cH_{\C^p}|_U \otimes \C[\xi_1 \dots \xi_q],
\qquad U \hbox{ open in } \C^p
$$
where $\C[\xi_1 \dots \xi_q]$ is the exterior algebra
generated by $\xi_1 \dots \xi_q$ and 
$\cH_{\C^p}$ denotes the sheaf of holomorphic functions on
$\C^p$.

A {\it complex supermanifold} of dimension $p|q$
is a superspace $M=(|M|, \cH_M)$
which is locally isomorphic to $\C^{p|q}$, i. e.
for all $x \in |M|$ there exist open sets $V_x \subset |M|$,
$U \subset \C^{p}$
such that:
$$
\cO_M|_{V_x} \cong \cH_{\C^{p|q}}|_U
$$ 

\end{definition}

\begin{definition}
A {\it superscheme} $S$ is a superspace $(|S|, \cO_S)$
such that $(|S|, \cO_{S,0})$ is a quasi-coherent sheaf of
$\cO_{S,1}$-modules. 
A {\it morphism} of supermanifolds  or of superschemes
is a morphisms of the corresponding superspaces.
\end{definition}

Superschemes can be characterized by a local model
as we shall presently see.

\begin{definition} $\uspec A$. \label{spec}

Let $A$ be a superalgebra and  
let $\cO_{A_0}$ be the structural sheaf of the ordinary scheme
$\uspec(A_0)=(\spec A_0, \cO_{A_0})$ ($\spec A_0$ denotes the
prime spectrum of the commutative ring $A_0$).
The stalk of the sheaf at the prime $p\in
\spec(A_0)$ is the localization of $A_0$ at $p$.
As for any superalgebra, $A$ is a
module over $A_0$. We have indeed a sheaf $\cO_A$ of
$\cO_{A_0}$-modules over $\spec A_0$ with stalk $A_{p}$, the localization
of the $A_0$-module $A$ over the prime
$p \in \spec(A_0)$:
$$
A_{p}=\{ {f \over g} \quad | \quad f\in A, g \in A_0-p\}.
$$
$A_p$ contains a unique two-sided maximal ideal 
generated by the maximal ideal in the local ring
$(A_p)_0$ and the generators of $(A_p)_1$ as $A_0$-module.

$\cO_A$ is a sheaf of superalgebras and
$(\spec A_0,\cO_A)$ is a superscheme that we 
denote with $\uspec A$. 
\end{definition}

The next proposition shows that $\uspec A$ is
the local model for superschemes.

\begin{proposition}
A superspace $S$ is a superscheme if and only if it is locally
isomorphic to $\uspec A$ for some superalgebra $A$, i. e.
for all $x \in |S|$, there exists $U_x \subset |S|$ open such
that $(U_x, \cO_S|_{U_x}) \cong \uspec A$. (Clearly $A$
depends on $U_x$).
\end{proposition}

\begin{proof} See \cite{cf} \S 3.
\end{proof}

\begin{definition}
We say that 
a superscheme $X$ is \textit{affine} if it is
isomorphic to $\uspec A$ for some algebra $A$
and we call $k[X]=_{\defi}A$ 
the \textit{coordinate ring} of the affine superscheme $X$. 
If $k[X]/I^{odd}$ is the coordinate ring of an ordinary affine algebraic
variety (called the \textit{reduced variety} 
associated to $X$) and $(|X|, \cO_{X,0})$ is a coherent sheaf of
$\cO_{X,1}$-modules, we say that $X$ is an {\it affine algebraic variety}.
\end{definition}


\begin{remark}
There is an equivalence of categories between superalgebras
and affine superschemes. This equivalence is treated in detail
in \cite{cf} \S 3.
\end{remark}

We now want to introduce the concept of functor of points
associated to an affine supervariety.

\begin{definition}
Let $X$ be a supervariety.
Its {\it functor of points} is given by:
$$
h_X: \salg \lra \sets, \qquad h_X(A)=\Hom (\uspec A,X)
$$
where $\salg$ is the category of commutative superalgebras.
If $X$ is an affine supervariety $h_X(A)=\Hom (k[X], A)$.
If $h_X$ is group valued we say that $X$ is an {\it affine supergroup}. 
This is equivalent to the fact that $k[X]$ is a Hopf superalgebra.
This is also the same as giving a multiplication
$m: X \times X \lra X$ and an inverse $i:X \lra X$ satisfying the
usual commutative diagrams.

\medskip

More in general, we say that $G: \salg \lra \sets$ is a
{\it supergroup functor} if it is group valued. Clearly, a representable
supergroup functor is an affine supergroup.
\end{definition}

\medskip

\section{Smoothness of complex algebraic supervarieties}

Let $k=\C$.

\medskip

Let $X=(|X|,\cO_X)$ be a supervariety
and let $P \in |X|$ be a \textit{closed point} i.e. $P$ corresponds
to a maximal ideal. Let $m_P$ be the maximal ideal in $\cO_{X,P}$.


\begin{definition} \label{smoothpt}
We say that $P$ is 
\textit{smooth} if 
$$
\widehat{\cO_{X,P}} \cong 
\C[[x_1 \dots x_r, \xi_1 \dots \xi_s]], \qquad
\widehat{\cO_{X,P}} =\lim_{\leftarrow} {\cO_{X,P}}/m_P^n 
$$
where $x_i$'s and $\xi_j$'s are respectively even and odd
variables.  
In this case we say that
the \textit{dimension} of the supervariety $X$ at $P$ is $r|s$.
Notice that the dimension is well defined, that is
if $\C[[x_1 \dots x_r, \xi_1 \dots \xi_s]]
\cong \C[[x_1 \dots x_m, \xi_1 \dots \xi_n]]$ then $r=m$, $n=s$.
\end{definition}

Smoothness of a point of a supervariety cannot be checked at the
classical level as the next examples show.

\begin{example} 

\medskip

\noindent
1. Consider the supervariety $X$ with  coordinate ring  
$\C[X]=\C[x,y,\xi,\eta]/(\xi\eta)$. Its reduced variety is the 
affine plane, where all the closed points are smooth in
the classical sense.
It is immediate to check that this supervariety has no smooth points
according to Definition \ref{smoothpt}.

\medskip

\noindent
2. Consider the supervariety with coordinate ring   
$\C[x,y,\xi,\eta]/(\xi x+\eta y)$. Again its reduced variety 
is the affine plane. One can check that all (closed) points are smooth 
except the origin.
\end{example}

Since the notion of smoothness is local we can assume that
$X$ is an affine supervariety, with coordinate ring 
$\C[X]=\C[x_1 \dots x_m,\xi_1 \dots \xi_n]/I$, where 
$I=(f_1 \dots f_p, \phi_1 \dots \phi_q)$.
In this case $\cO_{X,P}$ is the localization
of $\C[X]$ at the point $P$ (see Definition \ref{spec}). 

\medskip

\begin{definition}
As in the classical setting we define the 
\textit{jacobian} of  \break
$f_1 \dots f_p, \phi_1 \dots \phi_q$ 
$\in \C[x_1 \dots x_m,\xi_1 \dots \xi_n]$ at a point $P$ as:
$$
\Jac(f,\phi)=
\begin{pmatrix}{\del f_1  \over \del x_1}(P) & \dots & 
{\del f_1 \over x_m}(P)&
 {\del f_1  \over \del \xi_1}(P) & \dots & {\del f_1 \over \xi_n}(P) \\
\vdots & & \vdots & \vdots & & \vdots \\
{\del \phi_q  \over \del x_1}(P) & \dots & {\del \phi_q \over x_m}(P)&
 {\del \phi_q  \over \del \xi_1}(P) & \dots & {\del \phi_q \over \xi_n}(P) \\
\end{pmatrix}
$$
(for the definition of $\del f \over \del x$ see for example \cite{vsv}).
The \textit{rank} of the jacobian is given by $a|b$ where
$a$ and $b$ are the ranks of the $m \times p$, $n \times q$ diagonal
blocks.
\end{definition}

\begin{lemma} \label{jacobian}
Let the notation be as above. Let $P \in |X|$ be a closed point
i. e. a maximal ideal $m_P$ in $\C[X]$.
Then
$$
\rk(\Jac(f, \phi))(P)=
m|n-\Dim(m_P/m_P^2).
$$
\end{lemma}

\begin{proof}
The proof is the same as in ordinary case, (see for
example \cite{ha} pg 32), let's sketch it.
We have a natural identification:
$$
\begin{array}{ccc}
F: M_P / M_P^2 & \cong & \C^{m|n} \\
f & \mapsto & 
df_P=_{\defi}({\del f \over \del x_1}(P), \dots, {\del f \over \del \xi_n}(P))
\end{array}
$$
where $M_P$ denotes the maximal ideal corresponding to the point $P$ 
in \\ $\C[x_1 \dots x_m,\xi_1 \dots \xi_n]$.
Viewing the rows of  $\Jac(f,\phi)$ as vectors in
$\C^{m|n}$ the above identification tells us immediately that
$$
\rk (\Jac(f,\phi))(P)={\dim (I+M_P^2) / M_P^2}
$$
where $I=(f_1 \dots f_p,\phi_1 \dots \phi_q)$. Since 
localizations commute with quotients we have that:
$$
m_{P}/m_{P}^2 \cong (M_P/I) / ((M_P^2 +I)/I) =M_P/(M_P^2+I).
$$
Hence we have:
$$
\rk(\Jac(f,\phi))(P)={\dim (I +M_P^2) / M_P^2}=
\dim M_P/M_P^2-\dim M_P /(M_P^2+I).
$$
\end{proof}

\begin{proposition} \label{smooth}
If $P$ is a smooth point of an affine supervariety
$X$ with dimension $r|s$ in $P$ then:

\noindent
1. $m_P/m_P^2$ has dimension $r|s$. 

\noindent
2. $Gr(\cO_{X,P})=\C[x_1 \dots x_r, \xi_1 \dots \xi_s]$
where $Gr(\cO_{X,P})=\oplus_i m_P^i/m_P^{i+1}$.

\noindent
3. $\rk(\Jac(f, \phi)(P)=m|n-r|s$.

\end{proposition}

\begin{proof} Parts (1) and (2) are immediate
by Lemma \ref{maxideals} in the Appendix, (3) is a consequence
of Lemma \ref{jacobian}. 
\end{proof}

\begin{remark}
1. The proof of this result resembles the one for the commutative
setting. One difference that may ingenerate confusion is the
following. When we are localizing $\C[X]$ to obtain $\cO_{X,P}$ we are
using a maximal ideal of the even part $\C[X]_0$ that is
$(x_1-a_1, \dots x_m-a_m, \xi_i\xi_j, \forall i>j)$, $a_i \in \C$.
On the other hand, when we are completing the local 
superalgebra $\cO_{X,P}$ we are
taking the inverse limit of the system $\cO_{X,P}/m_P^n$, where
$m_P$ is the maximal ideal of this superalgebra, hence it is a graded
object and it will necessarily contain all the odd
generators.

\medskip

2. If $P$ is smooth, $m_P/m_P^2$ is generated by $r|s$ elements,
hence by the super Nakayama's Lemma \ref{snak},
we have that $m_P$ is generated by $r|s$ elements. 
\end{remark}

\begin{observation} \label{complexstructure}
The affine supervariety $X$ is embedded in $\C^{m|n}$ via the
chosen explicit presentation of its coordinate ring $\C[X]$. 
Hence we can give to the set of closed points of $X$ a complex
topology inherited from this embedding. However this topology
is independent from the embedding; this is a classical fact, still
valid in this setting since it is a topological question.
We want to show that the closed points of the supervariety $X$
equipped with this complex topology, admit a unique 
\textit{supermanifold structure} in a suitable 
complex neighbourhood $U$ of
the smooth point $P$. 
In other words we want to show that:
$$
\cH_{\C^{m|n}}|_U=(\cH_{\C^{m|n}}/\cK)|_U \cong 
\cH_{\C^{m|n}}|_U \otimes \C[\xi_1 \dots \xi_s] \qquad (*)
$$
where $\cK$ is the ideal sheaf whose global sections are
generated in $\cH_{\C^{m|n}}$ by the ideal $I$ of the supervariety $X$.
The whole question in the super setting if to show the existence of a local 
splitting $(*)$. To settle this problem
our strategy is to use the implicit functions
theorem, which is still valid in this setting. Let's recall
the statement from \cite{le} pg 52.
\end{observation}





\begin{theorem} \label{implicit}
Let $M$ be a complex supermanifold, $P \in |U|$, where $U \subset M$,
is isomorphic to an open in $\C^{r|s}$.
Let $K$ be the ideal in $\cH_M(U)$ generated by 
$g_1 \dots g_p, \gamma_1 \dots \gamma_q$ vanishing  
\footnote{We say that $f \in \cH_M(U)$ \textit{vanishes at}
$P$ if it is zero under the morphism:
$$
\cH_M(U) \lra \cH_{M,P} \lra \cH_{M,P}/m_{h,P} \cong \C
$$
$m_{h,P}$ being the maximal ideal in  $\cH_{M,P}$.} 
at $P$
and with linearly independent differentials at $P$. Then there exists
a unique subsupermanifold: 
$$
N=(|N|, \cH_N), \qquad \cH_N=\cH_M|_U/\cK
$$
where $\cK$ is the sheaf of ideals with global sections $K$
and $|N|$ is the topological space whose existence is granted by
the classical result.
\end{theorem}

\begin{remark}
The key for the proof of this result is the fact that any
set of functions $g_1 \dots g_p, \gamma_1 \dots \gamma_q$ with linearly
independent differentials at $P$ can be completed to obtain a
set of local coordinates in a neighbourhood of $P$. More details
on this can be found in \cite{vsv} pg 148.
\end{remark}

This theorem allows us, in a special case, to obtain immediately the
result we are after.

\begin{corollary} \label{completeintersection}
Let $P \in |X|$ be a smooth point, and let $X$ have dimension
$r|s$ at $P$.
Let's assume that the ideal $I$ of the supervariety $X$ is given by 
$I=(f_1 \dots f_{m-r}, \phi_1 \dots \phi_{n-s})$ 
(in this case we say that $X$ is  a complete intersection). Then
in a neighbourhood of  $P$, $X$ admits
a complex supermanifold structure (in the sense of Observation
\ref{complexstructure}).
\end{corollary}

\begin{proof}
This is a direct application of the Theorem \ref{implicit}.
The super variety $X$ is defined in $\C^{m|n}$ by the
polynomials $f_1 \dots f_{m-r}, \phi_1 \dots
\phi_{n-s}$ with $\rk(\Jac(f_i, \phi_{j})(P)=m|n-r|s$
Consider the ideal $K$ generated in $\cH(\C^{m|n})$ by the
$f_i$'s and $\phi_j$'s. Then there exists a unique 
subsupermanifold
$N$ of $\C^{m|n}$ such that $\cH_N=(\cH_{\C^{m|n}}/\cK)|_U$
for a suitable neighbourhood $U$ of $P$, $\cK$ is the
ideal sheaf whose global sections are $K$. 
\end{proof}

In general the ideal $I$ of the supervariety $X$ is given by 
$(f_1 \dots f_p,\phi_1 \dots \phi_q)$
where $p|q > m|n-r|s$.
We want to show that, as it happens for
the classical setting, $X$ is locally a complete intersection, so
that we can conclude our discussion with the same reasoning as
in Corollary \ref{completeintersection}.
Let $P \in X$ be a smooth point and 
assume $f_1 \dots f_{m-r},\phi_1 \dots \phi_{n-s}$ are such that: 
$$
\rk(\Jac(f_1 \dots f_{m-r},\phi_1 \dots \phi_{n-s}))(P)=m|n-r|s.
$$
Let $X'$ be the variety corresponding to the ring: 
$$
\C[X']=
\C[x_1 \dots x_m,\xi_1 \dots \xi_n]/
(f_1 \dots f_{m-r}, \phi_1 \dots \phi_{n-s})
$$ 
and let $\cO_{X',P}$ denote
its local ring at the closed point $P$.
We are going to show the following:

\medskip\noindent
1. $P$ is a smooth point of $X'$.
Moreover $X'$ has the same
dimension of $X$ i. e. $\cO_{X',P}=\C[[x_1 \dots x_r, \xi_1 \dots \xi_s]]$.
This implies that
$X'$ is a complete intersection.

\medskip\noindent
2. $X$ and $X'$ are locally isomorphic, in other words
$\cO_{X,P} \cong \cO_{X',P}$. Since
this result is true for all the points in a neighbourhood
of $P$, we have that $\cO_X(U) \cong \cO_{X'}(U)$.
Hence we can apply the result \ref{implicit} to $X'$ to conclude
that $X$ admits a supermanifold structure near $P$. 

\medskip

\begin{lemma} \label{surjective} Let the notation be
as above. We have the following commutative diagram:
$$
\begin{array}{ccc}
\cO_{X',P} & \twoheadrightarrow & \cO_{X,P} \\ \\
\downarrow & & \downarrow \\ \\
\widehat{\cO_{X',P}} & \twoheadrightarrow & \widehat{\cO_{X,P}}
\end{array}
$$
where the orizontal arrows are surjections, while
the vertical ones injections.
\end{lemma}


\begin{proof}
Observe that since we have a surjection $\C[X'] \lra \C[X]$
we also have a surjective morphism (this is a property of
localizations):
$$
\cO_{X',P} \lra \cO_{X,P}
$$
mapping the maximal ideal onto the maximal ideal.
This will give raise to a surjective system:
$$
\cO_{X',P}/{m_P'}^n \lra \cO_{X,P}/m_P^n.
$$
where $m_P$ and $m_P'$ denote the maximal ideals in
$\cO_{X,P}$ and $\cO_{X',P}$.
Hence $\widehat{\cO_{X',P}} \lra \widehat{\cO_{X,P}}$ is a
surjective map. The vertical arrows are injections since
$\cap m^i_P=\cap m'^i_P=(0)$. This happens since this is
true in the ordinary case and since the odd variables disappear
for large $i$'s.
\end{proof}

\begin{remark} 
By Lemma \ref{jacobian} we get immediately that
$\Dim (m_P/m_P^2)=\Dim (m_P'/m_P'^2)$ and since
the point $P$ is smooth
$$
\Dim (m_P/m_P^2) =\Dim (m_P'/m_P'^2) =r|s.
$$
Hence by the super Nakayama's
Lemma \ref{snak} we have that both $m_P$ and $m_P'$ are generated by
$r|s$ elements.
\end{remark}

\begin{lemma} \label{powserquotient}
Let the notation be as above. 
$$
\widehat{\cO_{X',P}} \cong \C[[x_1 \dots x_r, \xi_1 \dots \xi_s]]/I
$$
for a suitable ideal $I$.
\end{lemma}

\begin{proof}
By the Theorem \ref{universal} in the Appendix, we have that
there exist a unique map:
$$
\C[[x_1 \dots x_r, \xi_1 \dots \xi_s]] \lra \widehat{\cO_{X',P}}
$$
sending $x_i$'s and $\xi_j$'s into $r|s$ generators of
the maximal ideal $m'_P$. So the map is surjective and we
obtain our result.
\end{proof}

\begin{proposition} \label{seriesiso}
Let the notation be as above.
$$
\widehat{\cO_{X',P}} \cong \widehat{\cO_{X,P}}
=\C[[x_1 \dots x_r,\xi_1 \dots \xi_s]].
$$
\end{proposition}

\begin{proof}
By the Lemma \ref{surjective} we have that:
$$
\widehat{\cO_{X',P}}/J \cong \widehat{\cO_{X,P}}=
\C[[x_1 \dots x_r,\xi_1 \dots \xi_s]]
$$
By the Theorem \ref{powerseries} in the Appendix, 
we get the result.
\end{proof}

We have proven the local isomorphism in the completions, 
now we turn our attention to the local rings.

\begin{lemma} \label{localiso}
Let the notation be as above.
$$
\cO_{X',P} \cong \cO_{X,P}.
$$
\end{lemma}

\begin{proof} By Proposition \ref{seriesiso} and Lemma \ref{surjective}.



\end{proof}

This concludes the proof of the following:

\begin{theorem}
Let $X$ be a complex algebraic supervariety, $P$ a smooth point of $X$.
Then, there exist a neighbourhood of $P$ where we can give
to $X$ a unique structure of a complex supermanifold.
\end{theorem}

\begin{proof} Assume without loss of generality
that $X$ is affine and has dimension $r|s$ at $P$.
Let 
$$
\C[X]=\C[x_1 \dots x_m,\xi_1 \dots \xi_n]/(f_1 \dots f_p,\phi_1 \dots
\phi_q)
$$ 
be the coordinate ring of $X$. Let $X'$ be the algebraic
supervariety defined by the coordinate ring:
$$
\C[X']=\C[x_1 \dots x_m,\xi_1 \dots \xi_n]/(f_1 \dots f_{m-r},\phi_1 \dots
\phi_{n-s})
$$
where $\rk(\Jac(f_1 \dots f_{m-r},\phi_1 \dots
\phi_{n-s}))=m|n-r|s$.
Then by Corollary \ref{completeintersection} the result holds for $X'$ 
and by Lemma \ref{localiso} $X$ and $X'$ are locally isomorphic.
\end{proof}

The next lemma will be crucial in the discussion of smoothness
of algebraic supergroups in Section \ref{supergroups}.

\begin{lemma} \label{associatedgraded}
Let the notation be as above.
Let $m_P$ be generated by $r|s$ elements. 
If 
$$
Gr(\cO_{X,P})=\C[x_1 \dots x_r, \xi_1 \dots \xi_s]
$$
then $P$ is smooth.
\end{lemma}

\begin{proof}
By the Theorem \ref{universal} in the Appendix, we have
that there exist a surjective map:
$$
\C[[x_1 \dots  x_r,\xi_1 \dots \xi_s]] \lra \widehat{\cO_{X,P}}
$$ 
sending $x_1 \dots x_r,\xi_1 \dots \xi_s$ 
into the generators of the maximal ideal
of $\widehat{\cO_{X,P}}$. Hence we have that $\widehat{\cO_{X,P}}=
\C[[x_1 \dots  x_r,\xi_1 \dots \xi_s]]/J$ for some
ideal $J$.
Since $Gr({\cO_{X,P}})=Gr(\widehat{\cO_{X,P}})$ by
Lemma \ref{maxideals} the
 result follows by Lemma \ref{gr} in the Appendix.
\end{proof}

\section{Smoothness of Supergroups} \label{supergroups}

In this section we want to show that affine algebraic supergroups 
are smooth, that is all closed points are smooth.
In other words we show that the set of closed points
of an affine supergroup has a supermanifold structure in
the sense of Observation \ref{complexstructure}, hence it
is a Lie supergroup.
We will do this by using an argument appearing in
the classical Cartier's theorem which states that
Hopf algebras over a field of characteristic zero are reduced.

\medskip

It is enough to prove that the
identity is a smooth point, since, because 
of the multiplication law, all closed points have the same 
local structure.

\medskip

Let $G$ be an affine algebraic supergroup, 
$\C[G]$ its Hopf superalgebra with comultiplication $\Delta$,
counit $\epsilon$ and antipode $S$. 
Let $m_{1}= ker \epsilon$ be the maximal
ideal of the identity element and let:
$$
m_1 / m_1^2 = \Span_\C \{t_1 \dots t_{r+s}\} 
$$
where $t_1 \dots t_r$ are even, $t_{r+1} \dots t_{r+s}$ are
odd.
\medskip

By an abuse of notation let 
$t_1 \dots t_{r+s}$, 
denote also the image of these elements modulo $m_1^N$.

\medskip

\begin{lemma} \label{monomials}
The monomials $t_1^{n_1} \dots t_{r+s}^{n_{r+s}}$,
$\sum_{i=1}^{r+s} n_i=N$
form a basis for the superspace $m_1^N/m_1^{N+1}$.
(Clearly $n_i=0,1$ if $i$ is the index of an odd element,
$i=r+1 \dots s$). 
\end{lemma}

\begin{proof} The proof is the same as in the classical case,
we include it here for completeness (for more details
see \cite{wa} pg. 86).
Let $t_1^* \dots t_{r+s}^*$ be the dual basis of
$t_1 \dots t_{r+s}$. Define the map:
$$
d_l: \C[G]=\C \oplus m_1 \lra m_1 / m_1^2 \lra \C
$$
as $d_l=t_l^* \cdot p$, 
$l=1 \dots r+s$, where $p:\C[G] \lra m_1/m_1^2$ is the natural
projection.

Each $d_l$ gives rise to a derivation $D_l:\C[G] \lra \C[G]$
in the following way:
$$
D_l(a)=_{def}\sum a^{(1)}d_l(a^{(2)}), \qquad 
\hbox{where} \quad \Delta(a)=\sum a^{(1)} \otimes a^{(2)}.
$$
Observe that
$$
\ep (D_l(a))=\sum \ep(a^{(1)})d_l(a^{(2)})=
d_l \sum \ep(a^{(1)})a^{(2)}=d_l(a).
$$
Hence we have that $D_i(t_j) \equiv \delta_{ij}$ and 
modulo $m_1$ (since
$ker \ep = m_1$). 
Let $P(T_1 \dots T_{r+s})$ be a homogeneous polynomial 
of degree $n$ over $\C$. 
$$
D_i(P)(t_1 \dots t_{r+s})= 
\sum  D_i(t_j)
{\partial P \over \partial T_j}(t_1 \dots t_{r+s})
$$
Since ${\partial P \over \partial T_j}(t_1 \dots t_{r+s}) \in m_1^{n-1}$
we have $D_i(P) \equiv {\partial P \over \partial t_j}$ 
modulo $m_1^n$.
Now, since if $x \equiv y$ modulo $m_1^n$ it implies
$D_i(x) \equiv D_i(y)$ modulo $m_1^{n-1}$, we have that:
$$
\begin{array}{c}
D_{r+s}^{n_{r+s}} \dots D_1^{n_1}
t_1^{n_1} \dots t_{r+s}^{n_{r+s}}=
n_1! \dots n_{r+s}! \qquad \hbox{mod } m_1
\end{array}
$$
while on all other monomials the composition of $D$'s will give zero.
Hence given a relation $P$ in $m_1^{n+1}$ applying the correct
sequence of $D_i$'s one can single out the coefficient of any
monomial.
\end{proof}
\medskip

\begin{corollary}
The identity $1 \in |G|$ is a smooth point.
\end{corollary}

\begin{proof} 
By the  Lemma \ref{monomials} and Lemma \ref{maxideals}
we have that the graded associated
ring to $\cO_{G,1}$ is 
$$
Gr(\cO_{1,G})=\C[t_1 \dots t_r,\theta_1 \dots
\theta_s]. 
$$
This implies by Lemma \ref{associatedgraded} that 
the identity point is smooth. 
\end{proof}

\medskip 

\begin{corollary}
If $G$ is an affine supergroup, then $G$ is smooth, that is
all its closed points are smooth.
\end{corollary}

\begin{proof} Let $h_G$ denote the functor of points of $G$
and $\mu: h_G \times h_G \lra h_G$ the natural transformation
corresponding to the group law. Let $g \in |G|$ be a closed point. 
$g$ can be identified with an element of $h_G(\C) \subset h_G(A)$.
Hence we can define a natural transformation:
$$
l_g: h_G \lra h_G, \qquad l_{g,A}(x)=m_A(g,x),
\forall x \in h_G(A)
$$ 
This natural transformation corresponds to an isomorphism
of $G$ into itself, hence $\cO_{G,1} \cong \cO_{G,g}$,
so $g$ is smooth. (For more details on the correspondence
between natural transformations between functor of points
and morphisms of the supervarieties see \cite{cf} Chapter 3).
\end{proof}

\section{The Stabilizer Theorem}

Notation: In this section we use the same letter $X$ to denote
both a supervariety $X$ and its functor of points $h_X$.

\medskip

Let $G$ be an affine algebraic supergroup acting on an affine 
supervariety $X$,
in other words we have a morphism
$$
\rho:G \times X \lra X, \qquad (g,x) \mapsto  g \cdot x,
\qquad \forall g \in G(A), x \in X(A)
$$
satisfying the usual properties, viewed in the category of
supervarieties.
Let $u$ be a topological point of $X$, that is $u \in |X|$ or
equivalently $u \in X(\C)$ $=\Hom(\C[X], \C)$. Let $m_u$ be
the maximal ideal corresponding to $u$.
Notice that $u$ can be viewed
naturally as an $A$-point $u_A$ for all superalgebras $A$ since
$\C \subset A$. 
So we have a morphism:
$$
\tau:G \lra X, \qquad g \mapsto g \cdot u_A
$$
or equivalently:
$$
\ttau:\C[X] \lra \C[G].
$$

\begin{definition}
We define the {\it stabilizer supergroup functor} of the 
point $u \in |X|$ with respect to the action $\rho$, the 
group valued functor
$Stab_{u}: \salg \lra \sets$ defined by:
$$
Stab_{u}(A)=\{g \in G(A) \quad | \quad \tau_A(g)=g \cdot u_A = u_A\}
$$ 
where $\tau_A:G(A) \lra X(A)$,
or equivalently:
$$
\begin{array}{c}
Stab_{u}(A)=\{g \in G(A)=\Hom(\C[G],A) \quad | 
\quad g \cdot \ttau=u_A\} 
\end{array}
$$
\end{definition}

We want to prove that this functor is representable by an affine supergroup.



\begin{theorem} \label{stabthm}
Let $G$ be an affine supergroup acting on an affine supervariety $X$
and let $u$ be a topological point of $X$. Then $Stab_{u}$ is
an affine supergroup.
\end{theorem}

\begin{proof}
The stabilizer can be described in an equivalent way as:
$$
Stab_u(A)=\{g \in G(A) \quad | \quad (g \cdot \tau)|_{m_u}=0\}
$$
where $m_u \subset \C[X]$ is the ideal of $u$.
Let $I$ be the ideal in $\C[G]$ generated by $\ttau(x)$ for
all $x \in m_u$. One can immediately check that $g \in G(A)=
\Hom(\C[G], A)$ is in $Stab_u(A)$ if and only if $g$ factors
via $\C[G]/I$, that is $g: \C[G] \lra \C[G]/I \lra A$.
So we have that $Stab_u(A)=\Hom(\C[G]/I, A)$. 
\end{proof} 

We want to describe some important applications of this result.

\medskip


\section{The classical series of Lie supergroups}

In \cite{ka} Kac proved a classification theorem
for simple Lie superalgebras. The description of the
supergroup functors, corresponding to the classical
super series of Lie superalgebras introduced by Kac, appeared
in \cite{dm} pg 70; however no representation
statement was proved there.

\medskip
  
In this section we want to describe the supergroup functors 
corresponding to the classical super series and to show
they are representable i. e. they are algebraic supergroups, hence Lie
supergroups by the results of Section 4. For the series
$A(m,n)$, $B(m,n)$, $C(n)$ and $D(m,n)$ this result was
proved in \cite{vsv} pg 289 with differential tecniques.

\medskip
One should also prove that the Lie superalgebras\footnote{
For the definition of Lie superalgebra of an algebraic supergroup
see \cite{cf} Ch. 5} of these
Lie supergroups coincide with the classical series mentioned
above; however this goes beyond the scope of this paper and
we leave it as an exercise to the reader.
 

\medskip

1. {\bf $A(n)$ series}. 
Define $\rGL_{m|n}(A)$ as the set of all invertible
morphisms $g:A^{m|n}\rightarrow A^{m|n}$.
This is equivalent to ask 
that the {\it Berezinian} \cite{be} or {\it superdeterminant} 
$$
\Ber(g)=\Ber \begin{pmatrix} p & q \\ r & s \end{pmatrix} =
\det(p-qs^{-1}r)\det(s^{-1})
$$ 
is invertible in $A$ 
(where $p$ and $s$ are $m \times m$, $n \times n$ matrices
of even elements in $A$, while $q$ and $r$ are $m \times n$,
$n \times m$ matrices of odd elements in $A$). 
A necessary and sufficient 
condition for $g \in \rGL_{m|n}(A)$ to be invertible is that
$p$ and $s$ are invertible. The group valued functor
\begin{eqnarray*}
\rGL_{m|n}&:\salg &
\longrightarrow \sets\\ &A&\longrightarrow
\rGL_{m|n}(A).
\end{eqnarray*} 
is an affine supergroup called the \textit{general
linear supergroup} and it is represented by the algebra 
\begin{eqnarray*}&
\C[\rGL_{m|n}]:=\C[x_{ij}, y_{\alpha \beta},
\xi_{i\beta},\gamma_{\alpha
j},z,w]/\bigr((w\det(x)-1,z\det(y)-1\bigl),
\\ 
&i,j=1,\dots m,\;\;
\alpha,\beta=1,\dots n.
\end{eqnarray*} 

Consider the morphism:
$$
\rho: \rGL_{m|n} \times \C^{0|1} \lra \C^{0|1} \qquad
(g,c) \lra \Ber(g)c.
$$
The stabilizer of the point $1 \in \C^{0|1}$ 
coincides with all the matrices in $\rGL_{m|n}(A)$ with
Berezinian equal to 1, that is $\rSL_{m|n}(A)$ the
special linear supergroup. By the Theorem \ref{stabthm}
we have immediately that $\rSL_{m|n}$ is representable
and by the result of Section \ref{supergroups} we
have that it is a complex supermanifold.
Moreover one can check that $A(m,n)=\rLie(\rSL_{m|n})$.
\medskip

2. {\bf $B(m,n)$, $C(n)$, $D(m,n)$ series}.
Consider the morphism:
$$
\rho: \rGL_{m|2n} \times \cB \lra \cB \qquad
(g,\psi(\cdot,\cdot)) \lra \psi(g \cdot, g \cdot),
$$
where $\cB$ is the supervector space of all the
symmetric bilinear forms on $\C^{m|2n}$.
The stabilizer of the point $\Phi$ the standard
bilinear form on $\C^{m|2n}$ is the supergroup functor
$\rOsp_{m|2n}$. Again this is an algebraic supergroup
by Theorem \ref{stabthm} and it is also a complex
supermanifold. One can check that 
$B(m,n)=\rLie(\rOsp_{2m+1|2n})$, $C(n)=\rLie(\rOsp_{2|2n-2})$
and $D(m,n)=\rLie(\rOsp_{2m|2n})$.

\medskip

3. {\bf P(n) series}. 
Define the algebraic supergroup $\pi \rSp_{n|n}$ as we did
for $\rOsp_{m|n}$, by taking antisymmetric bilinear forms 
instead of symmetric ones. Consider the action:
$$
\begin{array}{ccc}
\pi \rSp_{n|n} \times \C^{1|0}  \lra  \C^{1|0} \qquad
(g,c)  \mapsto  \Ber(g)c.
\end{array}
$$
By Theorem \ref{stabthm} we have that $Stab_1$ is an affine
algebraic supergroup, hence it is a Lie supergroup. It is corresponding
to the $P(n)$ series.

\medskip

3. {\bf Q(n) series}. Let $D=\C[\eta]/(\eta^2+1)$. This is a
non commutative superalgebra. Define the supergroup functor
$GL_n(D):\salg \lra \sets$, with $GL_n(D)(A)$ the group  of automorphisms
of the left supermodule $A \otimes D$. In \cite{dm} is proven
the existence of a morphism called the \textit{odd determinant}
$$
\rodet: GL_n(D) \lra \C^{0|1}.
$$
Reasoning as before define:
$$
GL_n(D) \times \C^{0|1} \lra \C^{0|1}, \qquad g,c \lra \rodet(g)c. 
$$
Then $G=Stab_1$ is an affine algebraic supergroup and for $n \geq 2$
we define $Qg(n)$ as the quotient of $G$ and the diagonal subgroup
$GL_{1|0}$. This is an algebraic and Lie supergroup
and its Lie superalgebra is $Q(n)$.

\appendix
\section{Appendix: Commutative Superalgebra}

In this Appendix we collect some facts about commutative
superalgebra very similar to the
equivalent facts in commutative algebra. 

\medskip

Let $k$ be the ground field.

\medskip

All superalgebras are assumed to be commutative.
Let's denote (as before) with latin letter the even elements
and with greek letters the odd elements of a superalgebra.

\medskip

\begin{theorem}  \label{poly}
Let $A$ be a finitely generated superalgebra. 
Then there exists a unique superalgebra morphism 
$\phi: k[x_1 \dots x_m, \xi_1 \dots \xi_n] \lra A$ 
(where $k[x_1 \dots x_m, \xi_1 \dots \xi_n]$ denotes
the polynomial superalgebra with even indeterminates $x_i$'s
and odd indeterminates $\xi_j$'s)
sending the $x_i$'s and the $\xi_j$'s to chosen elements in $A$
of the correct parity.
\end{theorem}
This comes from the universality of the construction of the
polynomial superalgebra as it is done for example in
\cite{dm} pg 49.

\medskip

\begin{theorem}\label{universal}  
Let $A$ be a finitely generated superalgebra and let \break
$\hat A=\lim_{\longleftarrow}A/n^i$, be its completion with
respect an ideal $n$. Let $\hat n$ be the ideal in $\hat A$ 
corresponding to $n$.
Then there exist a unique superalgebra morphism
$\phi: k[[x_1 \dots x_m, \xi_1 \dots \xi_n]] \lra \hat A$ 
sending the $x_i$'s and the $\xi_j$'s to chosen elements in $\hat n$
of the correct parity.
\end{theorem}

\begin{proof} This is the same as Theorem 7.16 in \cite{ei} pg 200.
Let's briefly recall it.
By Theorem \ref{poly} we have that there is a unique map \break
$k[x_1 \dots x_m, \xi_1 \dots \xi_n] \lra \hat A/\hat n^i$ 
sending the $x_i$'s and the $\xi_j$'s to chosen elements in $n$.
Clearly this maps factors in the following way:
$$
k[x_1 \dots x_m, \xi_1 \dots \xi_n] \lra 
k[x_1 \dots x_m, \xi_1 \dots \xi_n]/(x_1 \dots x_m, \xi_1 \dots \xi_n)^i 
\lra \hat A/\hat n^i. 
$$
One can check that 
$$
{k[x_1 \dots x_m, \xi_1 \dots \xi_n]\over
(x_1 \dots x_m, \xi_1 \dots \xi_n)^i} 
\cong 
{k[[x_1 \dots x_m, \xi_1 \dots \xi_n]]
\over (x_1 \dots x_m, \xi_1 \dots \xi_n)^i} 
$$
hence by the universal property of the inverse limit we have obtained
the required map and the uniqueness.
\end{proof}


If $A$ is a local superring with maximal ideal $m$, 
let $Gr(A)=\oplus m^i/m^{i+1}$.

\begin{lemma} \label{gr}
If 
$$
Gr(k[[x_1 \dots x_r, \xi_1 \dots \xi_s]]/I) \cong 
Gr(k[[x_1 \dots x_r, \xi_1 \dots \xi_s]])
$$
then $I=(0)$.
\end{lemma}

\begin{proof} 
Let $m$ be the maximal ideal in $k[[x_1 \dots x_r, \xi_1 \dots \xi_s]]$.
There exist $i$ such that $I \subset m^i$ but $I \not\subset m^{i+1}$
otherwise we are done since $I \subset \cap m^i=(0)$.
Then
$$
(m^i/I)/(m^{i+1}+I/I)=m^i/(m^{i+1}+I) \neq m^i/m^{i+1}
$$
which gives a contradiction.
\end{proof}

\begin{theorem}\label{powerseries}
If 
$$
k[[x_1 \dots x_r, \xi_1 \dots \xi_s]]/I \cong 
k[[x_1 \dots x_r, \xi_1 \dots \xi_s]]
$$
then $I=(0)$.
\end{theorem}

\begin{proof} This is a consequence of Lemma \ref{gr}.
\end{proof}

\begin{lemma} \label{maxideals}
Let $A$ be a commutative superalgebra and $m$ a maximal ideal.
Let $A_m$ be the localization of $A$ into the even part
$m_0$ of the maximal ideal $m$ and $\widehat A_m$ the 
completion of $A_m$ with respect to the maximal ideal $\tilde m$
in $A_m$. Then:
$$
m^i /m^{i+1} \cong  \tilde m^i / \tilde m^{i+1} \cong 
\widehat m^i / \widehat m^{i+1}.
$$
\end{lemma}

\begin{proof} This is the same as in the commutative case,
because localization and completion commute with quotients.
\end{proof}

\begin{theorem} Super Nakayama's Lemma. \label{snak}
\medskip

Let $A$ be a local commutative super ring with maximal
(homogeneous) ideal $\mathfrak{m}$. Let $E$ be a finitely generated
module for the ungraded ring $A$.\\

\noindent (i) If $\mathfrak{m}E=E$, then $E=0$; more
generally, if $H$ is a submodule of $E$ such that $E=\frak mE+H$,
then $E=H$.\\

\noindent (ii) Let $(v_i)_{1\le i\le p}$ be a basis
for the $k$-vector space $E/\mathfrak{m}E$ where $k=A/\mathfrak{m}$. Let
$e_i\in E$ be above $v_i$. Then the $e_i$ generate $E$. If $E$ is
a supermodule for the super ring $A$, and $v_i$ are homogeneous
elements of the super vector space $E/\mathfrak{m}E$, we can choose the
$e_i$ to be homogeneous too (and hence of the same parity as
the $v_i$).\\

\noindent (iii) Suppose $E$ is projective, i.e. 
there is a $A$-module $F$ such that $E\oplus
F=A^N$ where $A^N$ is the free module for the ungraded ring $A$
of rank $N$. Then $E$ (and hence $F$) is free, and the
$e_i$ above form a basis for $E$.
\end{theorem}

\begin{proof} See \cite{cf} Appendix. \end{proof}

\end{document}